\magnification=1200
\input amssym.def
\input epsf
\def \BBZ {{\Bbb Z}}
\def \BBC {{\Bbb C}}
\def\sqr#1#2{{\vcenter{\vbox{\hrule height.#2pt
     \hbox{\vrule width.#2pt height#1pt \kern#1pt
     \vrule width.#2pt} \hrule height.#2pt}}}}
\def\square{\ ${\mathchoice\sqr34\sqr34\sqr{2.1}3\sqr{1.5}3}$}

\centerline {\bf ON KNOT INVARIANTS WHICH ARE NOT OF FINITE TYPE}
\bigskip
\centerline {
Theodore Stanford \footnote {*}
{Author partially supported by the
Naval Academy Research Council.}
and Rolland Trapp \footnote {**}
{Author partially supported by CSUSB Junior Faculty Support Grant.}}

\bigskip
\leftskip=1truein
\rightskip=1truein
\noindent
{\bf Abstract.} We observe that most known results of the form
``$v$ is not a finite-type invariant'' follow from two basic
theorems.  Among those invariants which are not of finite type,
we discuss examples which are ``ft-independent'' and examples
which are not.  We introduce $(n,q)$-finite invariants, which
are generalizations of finite-type invariants based on Fox's
$(n,q)$ congruence classes of knots.

\leftskip=0in
\rightskip=0in
\bigskip
\noindent
{\bf 1. How not to be of finite type.}
\medskip

Vassiliev~[21] defined a family of knot invariants based on
the study of singular knots.  Gusarov~[9] independently described
the same family of invariants using very different methods
and these invariants, now often called {\it finite-type
invariants}, have been derived and analyzed from a number of
different points of view.  There is now quite a large body
of literature.  The widespread interest is due mostly to
the following theorem, various parts of which were proved
by various authors.  See Birman~[3] for an exposition
and general proof.  See also Bar-Natan~[2], Birman and Lin~[4],
and Gusarov~[8].  By ``quantum knot polynomial''
we mean the Jones polynomial or one of its many generalizations.

\bigskip
\noindent
{\bf Theorem A.} \ Let $P_K (t)$ be a quantum knot polynomial,
and let $a_n (K)$ be the $x^n$ coefficient of the power
series obtained by substituting $t = e^x$ into $P_K(t)$.
Then $a_n$ is a finite-type invariant of order $\le n$.

\medskip
The invariants in Theorem~A will all be rational-valued.
However, we take ``finite-type invariant'' to mean any
knot invariant $v$ taking values in an abelian group $A$,
such that there exists a positive integer $n$ with
$v (K) = 0$ for any singular knot $K$ with more than $n$
singularities.  The value of $v$ on a singular knot is
determined by taking the difference between the positive
and negative resolutions of the singularity in the standard
way.

Theorem~A gives us many finite-type
invariants, and an immediate
and natural question becomes
whether there are invariants
which are not of finite type.
It is well-known that all finite-type invariants
are determined by those which take values in
$\BBZ$ and those which take values in $\BBZ_m$
(for all $m$), though in fact all known examples
are determined by the $\BBZ$-valued invariants.
The set of finite-type invariants taking values
in all these abelian groups is easily seen to
be countable, whereas the set of all $\BBZ$-valued
knot invariants is uncountable.
So in fact
there are many invariants which are not of finite type.
One could also ask whether there are invariants which are
not determined by finite-type invariants.  This is the
same question as whether there exists a pair of knots,
all of whose finite-type invariants are equal.  The answer
to this question is unknown.

Along similar lines, one may ask whether a particular
known and studied knot invariant is of finite type.
A number of standard knot invariants have been shown
not to be of finite type by different authors using
various techniques.
See Altschuler~[1],
Birman and Lin~[4], Dean~[5], Eisermann~[6], Ng~[13], and
Trapp~[20].
Most of these results follow from one or both of the
following two theorems:

\bigskip
\noindent
{\bf Theorem B.} \ No invariant taking a unique value on the
unknot is of finite type.  That is, if $v(K) = v ({\rm unknot})$
implies that $K$ is the unknot, then $v$ is not a finite-type
invariant.

\medskip
\noindent
{\it Proof:} \ There are now a number of constructions yielding
nontrivial knots with trivial invariants up to any fixed
order $n$.  See for example Gusarov~[9], Lin~[11],
Ohyama~[14], Stanford~[18].
\square

\medskip
It follows from Theorem~B that none of the following are
finite-type invariants:
crossing number, unknotting number, bridge number,
tunnel number, braid index, genus, free genus, Seifert genus,
stick number, and crookedness.

\bigskip
\noindent
{\bf Theorem C.} \ Any knot invariant $v$ for
which connected sums
cannot cancel is not of finite type.
That is, if $v$ is a knot invariant such that there exists
a knot $K$ with $v (K \# K^\prime) \ne v ({\rm unknot})$
for all knots $K^\prime$, then $v$ is not a finite-type
invariant.

\medskip
\noindent
{\it Proof:} \ Gusarov~[9] showed that for any knot $K$ and
any integer $n$, there exists a knot $K^\prime$ such that
$v (K\#K^\prime) = v ({\rm unknot})$ for any finite-type invariant
of order $<n$.  Other proofs of this result may be found in
Habiro~[10], Ng~[13], and Stanford~[19].
\square

\medskip
Let $P_K (t)$ be any knot polynomial with integer coefficients
which has the property that
$P_{K \# K^\prime} = C P_K (t) P_{K^\prime} (t)$, where
$C \in \BBZ[t^{\pm 1}]$
is some fixed constant.  (Most of the quantum knot polynomials
fit this description.)
Then by Theorem~C none of the following are finite-type
invariants: $P_K(t)$ as an element of the abelian group of
Laurent polynomials; $P_K (n)$, where $n$ is any integer
such that there exists a knot $K$ with $P_K (n) \ne \pm 1$;
$P_K (\alpha)$, for any non-algebraic $\alpha \in \BBC$; the
span of $P_K (t)$; the first and last coefficients of
$P_K(x)$; the degrees of the first and last coefficients of
$P_K (t)$.

Let $G$ be a nonabelian group, and let $v$ be a knot invariant
with some value $a$ such that $v(K)=a$ implies the existence
of an nonabelian homomorphism from the fundamental group
of the complement of $K$ into $G$.  By Theorem~C, $v$ is
not a finite-type invariant, because if such a homomorphism
exists for $K$ then it is easy to see that it exists for
$K \# K^\prime$ no matter what $K^\prime$ is.  If $G$ is
finite, then the number of homomorphisms from the group
of $K$ into $G$ is not a finite-type invariant. Nor is the
determinant of a knot (the Alexander polynomial evaluated at
$t=-1$), since if this is not $1$ then there exist
nontrivial representations of the group of $K$ into some
dihedral group.

\bigskip
\noindent
{\bf 2.  Ft independence.}
\medskip

A well-known invariant which is not covered by either
Theorem~B or Theorem~C is the signature of a knot.
The signature was originally shown not to be of finite type
by Dean~[5] and Trapp~[20].  Ng~[13] has proved
a much more general theorem:

\bigskip
\noindent
{\bf Theorem D.} \ The only knot concordance invariant
which is also of finite type is
the Arf invariant.  In fact, for any knots $K_1$ and $K_2$
with equal Arf invariant, and for any positive integer $n$,
there exists a knot $K_3$ with the same concordance class as
$K_1$, and such that $v (K_2) = v(K_3)$ for any finite-type
invariant of order $<n$.

\medskip
Theorem~D inspires a definition: \ we shall say that a knot
invariant $v$ is {\it ft-independent} if for any
set of finite-type invariants $v_1, v_2, \dots v_m$, the
values of $v_i(K)$ put no restrictions on the values of
$v(K)$.  More precisely, $v$ is ft-independent if
for any knots $K_1$ and $K_2$ and for any positive integer
$n$ there exists a knot $K_3$ such that $w (K_3) = w(K_2)$
for any finite-type invariant $w$ of order $<n$, and
such that $v (K_3) = v(K_1)$.

Thus, by Theorem~D, any knot concordance invariant which is
independent of the Arf invariant is ft-independent.  However,
many of the invariants that we have listed as being not of
finite-type as a result of Theorems~B and~C are also clearly
not ft-independent.  Crossing number, for example, is not
ft-independent.  Choose a positive integer $m$ and a
finite-type invariant $v$.  If $v(K) \ne v(K^\prime)$
for all knots $K^\prime$ of crossing number $\le m$, then
clearly the crossing number of $K$ is greater than $m$.
Thus the value of $v$ places some restrictions on the possible
crossing number of a knot.

Another invariant which is neither finite-type nor
ft-independent is the degree of the Conway polynomial.  If
the $(2n)$th coefficient is nonzero then the degree of the
Conway polynomial is at least $2n$.  The coefficients of the
Conway polynomial are themselves finite-type invariants
(Bar-Natan~[2]), the $t = e^x$ substitution being
unnecessary because the Conway polynomial doesn't have terms
of negative degree.  Hence nonzero values for these
invariants place lower bounds on the degree of the Conway
polynomial.  Since the degree of the Conway polynomial
(divided by $2$) is a lower bound for the genus of a knot,
it follows that the genus is also not ft-independent.

Clearly, if $P_K (t)$ is a quantum knot polynomial, then,
as an invariant taking values in the abelian group 
of Laurent polynomials
with integer coefficients, it is not ft-independent because
of Theorem~A.

One example of an ft-independent invariant which is not
a concordance invariant is the number of prime factors of
a knot. Given a knot $K$ and
a positive integer $n$, one can form the connected sum
of $K$ with any
number of knots whose invariants are trivial up to order $n$
so as to produce a $K^\prime$ with an arbitrarily large
number of prime factors.  Going the other way, it was
shown in Stanford~[18] that given a knot $K$ and a positive
integer $n$, there exists a prime knot $K^\prime$ such
that $v (K) = v (K^\prime)$ for every finite-type invariant
of order $\le n$.

\bigskip
\noindent
{\bf 3.  $(n,q)$-finite invariants.}
\bigskip

In this section we use Fox's notion of $(n,q)$ congruence of
knots to generalize finite-type invariants.  Before defining
$(n,q)$-finite invariants, however, we recall the definition of
$(n,q)$-congruence classes of links (see Fox [7], Nakanishi and
Suzuki [12], Przytycki [15]).  
We consider only links in ${\bf S}^3$,
and follow the definition given in Przytycki [16].

Given a link $L$ and a disk $D^2$ which $L$ intersects
transversely, let $U = \partial D^2$ and $q = | lk(L,U) |$.
Note that the complement of $U$ in ${\bf S}^3$ is a solid
torus $T$ with meridinal disk $D^2$.  A $t_{2,q}$ move on
$L$ is the restriction to $L$ of a Dehn twist in $T$ on the
disk $D^2$.  Thus a $t_{2,q}$ move has the effect of cutting
$L$ along $D^2$, inserting a full twist and reglueing.  A
$t_{2n,q}$ move is just the result of $n$ Dehn twists on
$D^2$, i.e. $t_{2n,q} = t^n_{2,q}$.  Two links $L_1,\ L_2$
are then called 
{\it congruent modulo $(n,q)$}, denoted $L_1 \equiv L_2
({\rm mod}\ n,q)$, if one can obtain $L_2$ from $L_1$ by a
sequence of $t_{2n,q'}^{\pm 1}$ moves together with
isotopies, where the $q'$ can vary but are required to be
multiples of $q$.

The Alexander Module is a good tool for constructing
invariants of $(n,q)$ congruence, and this is done by Fox as
well as Nakanishi and Suzuki.  In Przytycki [16], the
$t_{2n,q}$ moves are considered as generalizations of crossing
changes.  Generalized unknotting numbers are defined, and
lower bounds for these unknotting numbers are given.
Analogously, we will define generalized finite-type
invariants by replacing the crossing changes in the
Vassiliev skein relation by the more general
$t_{2n,q}$ moves.

Recall
the diagramatic definition of finite-type invariants.
A link invariant $f$ is of order $\le m$ 
if for each link diagram $D$
and every collection $(a_1,\dots,a_{m +1})$ of 
$m +1$ crossings, the alternating sum
$$\sum_{\vec{i}}(-1)^{|\vec{i}|}f(D_{\vec{i}}) \eqno {(3.1)}$$ 
vanishes,
where $\vec{i}$ is in $\BBZ_2^{m+1}$, $|\vec{i}|$ is the
number of nonzero coordinates in $\vec{i}$, and
$D_{\vec{i}}$ is the diagram $D$ with crossing $a_j$ changed
whenever the $j^{th}$ coordinate of $\vec{i}$ is one.
Thinking of $t_{2n,q}$ moves as generalized crossing
changes, we define $f$ to be 
{\it $(n,q)$-finite of order $\le m$} if
for all links $L$ and collections $D^2_1,\dots,D^2_{m+1}$ of
mutually disjoint disks, with $\partial D^2_i = U_i$,
satisfying $q_j = {\rm lk}(U_j, L) \equiv 0\ {\rm mod}\ q$,
the alternating sum
$$\sum_{\vec{i}}(-1)^{|\vec{i}|}f(L_{\vec{i}})\eqno{(3.2)}$$
vanishes, where $L_{\vec{i}}$ is obtained by
$t_{2n,q_j}$ moves on each $U_j$ for which the $j^{th}$
coordinate of $\vec{i}$ is one.

\bigskip
\noindent
{\bf Theorem E.} \ Let $q \ge 0$.  If $f$ is a finite-type
invariant of order $\le m$, then $f$ is an $(n,q)$-finite
invariant of order $\le m$.  If $f$ is $(1,q)$-finite
of order $\le m$, then $f$ is finite-type of order $\le m$.

\medskip
\noindent
{\it Proof.} \ Let $f$ be a $(1,q)$-finite invariant of order
$\le m$.  In our definition of $(n,q)$-finite, the $q_j$ are
allowed to be $0$.  Since crossing changes can be obtained
by a Dehn twist on a disk whose boundary has linking number
zero with the link, any sum (3.1) can be written in the form
of (3.2) (with $n=1$).  Thus $f$ is of finite type.

Now suppose that $f$ is a finite-type invariant of order
$\le m$.  Note that the difference $f(L) - f(t_{2,q}(L))$ can be
realized as a sum of values of $f$ on singular links.  This
follows directly from Theorem 3.1 of Stanford [17], but the
intuitive idea is that the full twist in $t_{2,q}(L)$ can be
undone by crossing changes, and thus any sum (3.2)
may be written as a sum of expressions (3.1).  This
argument works for $t_{2n,q}$ as well, or else one can
note that an $(n,q)$-finite invariant 
is $(kn,q)$-finite for all $k \in \BBZ^+$, because 
$t_{2kn,q} = t_{2n,q}^k$.  \square
\medskip

Thus $(n,q)$-finite invariants include finite-type
invariants.  It remains to show that the inclusion is
proper.  It is easy to see that invariants of
$(n,q)$-congruence classes of links are exactly the
$(n,q)$-finite invariants of order $0$ (in the same
way that the only finite-type invariants of order $0$
are those which are invariant under crossing changes,
namely, anything which depends only on the number of
components in the link).  Thus we need only find an
invariant of $(n,q)$-congruence which is not of
finite type.  The following is a direct result of
Lemma 2.6b of Przytycki [16], which states that the number
of $2n$-colorings of a link is invariant under 
certain $t_{2n,q}$ moves, and of Theorem~C above,
which implies that the number of colorings of a knot
is not a finite-type invariant:

\bigskip
\noindent
{\bf Theorem F.} \ The number of $2n$ colorings of a link
is an order $0$ $(n,q)$-finite invariant for any even $q$.
Thus there are order $0$ $(n,q)$-finite invariants which are
not of finite type.

\bigskip
\noindent
{\bf Remark:}  Nakanishi and Suzuki define two links to be 
$q-congruent\ modulo\ n$ if they are related by a sequence of
$t_{2n,q}^{\pm 1}$ moves, thus making the restriction that 
$lk(U,L) = q$.  One can similarly refine the
definition of $(n,q)$-finite invariants.  Lemma 2.6a of 
Przytycki [15]
can be applied to this refined notion of finite invariants.  The
reason it doesn't apply directly to $(n,q)$-finite invariants as
defined above is that one could have $lk(U,L) 
= 0$.  This implies that the geometric linking number is even
and violates the hypotheses of Lemma 2.6a in [15].

\bigskip
\noindent
{\bf References.}
\medskip

\smallskip
\item {[1]}
Daniel Altschuler,
{\it Representations of knot groups and Vassiliev invariants},
Journal of Knot Theory and its Ramifications
5 (1996), no. 4, 421--425.

\smallskip
\item{[2]}
D. Bar-Natan, {\it On the Vassiliev knot invariants},
Topology  34 (1995) no. 2, 423--472.

\smallskip
\item {[3]}
J. S. Birman, {\it New points of view in knot theory},
Bulletin of the American Mathematical Society
28 (1993) no. 2, 253--287.

\smallskip
\item{[4]}
J.S. Birman and X.-S. Lin, {\it Knot polynomials and
Vassiliev's invariants}, Inventiones mathematicae 111 (1993) no. 2,
225--270.

\smallskip
\item {[5]}
J. Dean,
{\it Many classical knot invariants are not Vassiliev invariants},
Journal of Knot Theory and its Ramifications 3 (1994) no. 1, 7--10.

\smallskip
\item {[6]}
Michael Eisermann,
{\it The number of knot group
representations is not a Vassiliev invariant},
Preprint.

\smallskip
\item{[7]}
R. H. Fox,
{\it Congruence Classes of Knots},
Osaka Math. J. 10 (1958), 37--41.

\smallskip
\item {[8]}
M.N. Gusarov,
{\it A new form of the Conway-Jones polynomial of oriented links},
Topology of manifolds and varieties, 167--172.
Advances in Soviet Mathematics 18,
American Mathematical Society, 1994.

\smallskip
\item {[9]}
M.N. Gusarov,
{\it On $n$-equivalence of knots and invariants of finite degree},
Topology of manifolds and varieties, 173--192,
Advances in Soviet Mathematics 18,
American Mathematical Society, 1994.

\smallskip
\item {[10]}
K. Habiro, {\it Claspers and the Vassiliev skein modules},
preprint, University of Tokyo.

\smallskip
\item {[11]}
X.-S. Lin,
{\it Finite type link invariants of $3$-manifolds},
Topology 33 (1994), no. 1, 45--71.

\smallskip
\item{[12]}
Y. Nakanishi and S. Suzuki,
{\it On Fox's Congruence Classes of Knots},
Osaka J. Math 24 (1987), 217--225.

\smallskip
\item {[13]}
K.Y. Ng,
{\it Groups of ribbon knots},
Topology 37 (1998), no. 2, 441--458.

\smallskip
\item {[14]}
Y. Ohyama,
{\it Vassiliev invariants and similarity of knots},
Proceedings of the American Mathematical Society
123 (1995) no. 1, 287--291.

\smallskip
\item{[15]}
J. Przytycki,
{\it 3-Coloring and Other Elementary Invariants of Knots}.

\smallskip
\item{[16]}
J. Przytycki, 
{\it $t_k$-moves on links},
Contemporary Mathematics, Volume 78 (1988), 615--655.

\smallskip
\item{[17]}
T. Stanford,
{\it The Functorality of Vassiliev-Type Invariants of Links,
Braids, and Knotted Graphs},
Random knotting and linking (Vancouver, BC, 1993).
Journal of Knot Theory and its Ramifications 3 (1994),
no. 3, 15--30.

\smallskip
\item {[18]}
T. Stanford,
{\it Braid commutators and Vassiliev invariants},
Pacific Journal of Mathematics 174 (1996)
no. 1, 269--276.

\smallskip
\item {[19]}
T.B. Stanford,
{\it Vassiliev invariants and knots modulo pure braid subgroups},
Preprint {\tt GT/9805092} available from
{\tt front.math.ucdavis.edu}.

\smallskip
\item {[20]}
R. Trapp,
{\it Twist sequences and Vassiliev invariants},
Random knotting and linking (Vancouver, BC, 1993).
Journal of Knot Theory and its Ramifications 3 (1994),
no. 3, 391--405.

\smallskip
\item{[21]}
V. A. Vassiliev, {\it Cohomology of knot spaces},
Theory of Singularities and Its Applications, 23--69,
Advances in Soviet Mathematics 1,
American Mathematical Society, 1990.

\bigskip
\noindent
Mathematics Department

\noindent
United States Naval Academy

\noindent
572 Holloway Road

\noindent
Annapolis, MD \ 21402

\smallskip
\noindent
{\tt stanford@nadn.navy.mil}

\bigskip
\noindent
Mathematics Department

\noindent
California State University, San Bernardino

\noindent
San Bernardino, CA \ 92407

\smallskip
\noindent
{\tt trapp@math.csusb.edu}

\end